\documentclass[a4paper,11pt,oneside]{article}

\usepackage{amsmath,amsthm,amsfonts,amssymb,amscd}
\usepackage{mathtools}
\usepackage{pgf,tikz}
\usepackage{float}
\usepackage{enumerate}
\usepackage{booktabs}
\usepackage{pgfplots}
\usepackage[noend]{algpseudocode}
\usepackage{graphicx,bm,xcolor}
\usepackage{algorithm}
\usepackage{algpseudocode}
\usepackage{hyperref}
\usepackage{caption}
\usepackage[T1]{fontenc}
\usepackage{t1enc}
\usepackage[polish,german,english]{babel}
\usepackage{verbatim} 

\usepackage{url}
\usepackage{authblk}

\usepackage{cleveref}
\usepackage{ao-math-std}
\usepackage{ao-math-fields}
\usepackage{newunicodechar}

\newcommand\Vu{V_u}
\newcommand\Vphi{V_\varphi}
\newcommand\Yu{Y_u}
\newcommand\Yphi{Y_\varphi}

\newcommand\aein{\ifmmode\text{ a.e.\ in }\else a.e.\ in \fi}
\newcommand\aetin{\ifmmode\text{ for a.e.\ $t$ in }\else a.e.\ t\in \fi}
\newcommand\pd{\partial}
\newcommand\bl{\bar l}
\newcommand\bq{\bar q}
\newcommand\bu{\bar u}
\newcommand\by{\bar y}
\newcommand\bphi{\bar\varphi}
\newcommand\bbl{\bar\boldl}
\newcommand\bbu{\bar\boldu}
\newcommand\setE{\mathcal E}
\newcommand\setI{\mathcal I}
\newcommand\mcG{\mathcal{G}}
\newcommand\mcY{\mathcal{Y}}
\newcommand\mcZ{\mathcal{Z}}
\newcommand\mcK{\mathcal{K}}

\newcommand\wu{w_u}
\newcommand\wphi{w_\varphi}

\newcommand\qfor{\qtextq{for all}}
\newcommand\tfor{\text{ for all }}
\let\bs\boldsymbol
\newcommand\boldl{\bs l}
\newcommand\boldu{\bs u}

\newcommand\boldz{\bs z}

\newcommand\boldPhi{\bs\Phi}
\newcommand\bolddeltal{\bs{\delta l}}
\newcommand\bolddeltau{\bs{\delta u}}

\newcommand\Phiu[1][]{\Phi_{u\optindex{#1}}}

\newcommand\Phiphi[1][]{\Phi_{\varphi\optindex{#1}}}

\newcommand\dPhiphi{\dot{\Phi}_\varphi}
\newcommand\gkappa{\textsl{g}_{\kappa}}

\newcommand\J{\mathcal J}
\newcommand\Lagr{\mathcal L}
\newcommand\Linf{L^\infty}
\newcommand\dphi{\dot\varphi}

\newcommand\CIQ{W}

\newcommand\diff{\,\mathrm{d}}
\newcommand\ds{\diff s}
\newcommand\dt{\diff t}
\newcommand\dx{\diff x}

\newcommand\IO{I \times \Omega}
\newcommand\iprod[3][]{(\fcarg{#2},\fcarg{#3})\optsb{#1}}
\newcommand\dualprod{\sprod[\Vphi^*,\Vphi]}
\newcommand\bidualprod{\sprod[\Vphi^{**},\Vphi^*]}

\newcommand\Idualprod[2]{\int_I \dualprod{#1}{#2} \dt}

\newcommand\prodO{\iprod[\Omega]}
\newcommand\prodIN{\iprod[I \times \Gamma_N]}
\newcommand\prodIO{\iprod[\IO]}

\newcommand\normIO{\norm[\IO]}

\DeclareMathOperator\cone{cone}

\DeclareMathOperator\tint{int}
\DeclareMathOperator\Span{span}

\newunicodechar{φ}{\varphi}
\newunicodechar{Δ}{\Delta}
\newunicodechar{∂}{\partial}
\newunicodechar{χ}{\chi}
\newunicodechar{ψ}{\psi}
\newunicodechar{Φ}{\Phi}
\newunicodechar{κ}{\kappa}
\newunicodechar{ε}{\varepsilon}
\newunicodechar{μ}{\mu}
\newunicodechar{λ}{\lambda}
\newunicodechar{ℝ}{\mathbb{R}}
\newunicodechar{δ}{\delta}
\newunicodechar{∇}{\nabla}
\newunicodechar{ℕ}{\mathbb{N}}
\newunicodechar{τ}{\tau}
\newunicodechar{η}{\eta}
\newcommand{\Ombar}{\overline{\Omega}}

\theoremstyle{plain}
\newtheorem{theorem}{Theorem}[section]
\newtheorem{proposition}{Proposition}[section]
\newtheorem{lemma}[theorem]{Lemma}

\newtheorem{Remark}[theorem]{Remark}
\newtheorem{definition}[theorem]{Definition}

\newtheorem{Problem}[theorem]{Problem}

\definecolor{LUH-red}{RGB}{229,0,27}
\definecolor{LUH-orange}{RGB}{237,103,1} 

\definecolor{gruen}{RGB}{10,190,10}
\definecolor{purple(x11)}{rgb}{0.63, 0.36, 0.94}

\usepackage[top=2cm, bottom=2cm, left=2cm, right=2cm]{geometry}

\begin{document}

\title{Analysis of a space-time phase-field fracture complementarity
  model and its optimal control formulation}
\author[1]{D. Khimin}
\author[1]{J. Lankeit}
\author[1]{M.C. Steinbach}
\author[1]{T. Wick}

\affil[1]{Leibniz Universit\"at Hannover, Institut f\"ur Angewandte
  Mathematik, Welfengarten 1, 30167 Hannover, Germany}

\date{}

\maketitle
	
\begin{abstract}
The purpose of this work is the formulation of optimality conditions for phase-field optimal
control problems. The forward problem is first stated as an abstract nonlinear optimization problem,
and then the necessary optimality conditions are derived. The sufficient optimality conditions are also examined.
The choice of suitable function spaces to ensure the regularity of the nonlinear
optimization problem is a true challenge here.
Afterwards the optimal control problem with a tracking type cost functional is formulated.
The constraints are given by the previously derived first order optimality conditions of the forward problem.
Herein regularity is proven under certain conditions and first order optimality conditions
are formulated.
\textbf{AMS:} 49J50, 49K20, 74R10, 49J40
\end{abstract}

\section{Introduction}
Variational phase-field methods for the modeling of
fracture propagation are an important research area in applied mathematics and engineering.
First works establishing phase-field methods for fracture propagation
from a mathematical and mechanical point of view are
\cite{BourFraMar00,Bour07,BuOrSue10,MieWelHof10a,KuMue10}.
There are numerous additional references cited in the overview articles and monographs
\cite{BourFraMar08,BouFra19,AmGeraLoren15,WuNgNgSuBoSi19,Wi20_book,Fra21,DiLiWiTy22}
as well. The large majority of studies is concerned with forward modeling of phase-field fracture
with applications in numerous fields.
Since the year $2017$, optimization
with phase-field fracture as forward problem is being investigated
within optimal control
\cite{HehlNeitzel:2022,MoWo21,NeiWiWo17,NeiWiWo19,KhiSteiWi22_JCP,Hehl2022,KhiSteiWi23_JOTA}
as well as topology optimization \cite{DESAI2022111048}
and stochastic phase-field fracture settings \cite{GERASIMOV2020113353}.

Phase-field fracture forward problems can be classified as coupled
variational inequality systems (CVIS) \cite{Wi20_book} in which vector-valued
displacements couple with a smoothed indicator phase-field function.
In situations where fracture healing
is not allowed, as in our work and most often the case in the literature,
the phase-field function is subject
to an inequality constraint in time. Various approaches have been employed
to represent the inequality constraint, such as imposing Dirichlet values
in the fracture zone \cite{BourFraMar00,Bour07}, strain history function \cite{MieWelHof10b},
simple penalization \cite{MiWheWi15b,GERASIMOV2019990}, augmented Lagrangian formulations
\cite{WheWiWo14}, a closely related inexact augmented Lagrangian method \cite{Wi17_SISC},
primal-dual active set methods \cite{HeWheWi15},
interior-point methods
\cite{WAMBACQ2021113612}, recursive multilevel trust region methods in which
the corrections satisfy the irreversibility condition \cite{KoKr20}, truncated nonsmooth
Newton multigrid methods in which the variational structure handles the
irreversibility constraint pointwise \cite{GrKieSa23},
and complementarity formulations \cite{MaWiWo20}. In this work,
the latter is of interest for which we notice that \cite{MaWiWo20} formally
introduced and implemented a complementarity condition,
but without rigorous mathematical analysis and still in time-incremental form, thus
not within a space-time setting. In this respect, we notice that the first
space-time phase-field fracture formulation (with penalization of the irreversibility constraint)
as forward problem and within optimal control was proposed in \cite{KhiSteiWi22_JCP}.

The main objective of this work is the rigorous investigation of optimality conditions in terms of
KKT (Karush-Kuhn-Tucker) systems for phase-field fracture forward and optimal control problems
in a continuous space-time setting.
In optimal control, a cost functional
is minimized subject to some forward problem
that acts as a (physical) constraint.
Within this upper level control problem,
the fracture is driven by the
control, which can act as a boundary condition or a right hand side force
\cite{HehlNeitzel:2022,MoWo21,NeiWiWo17,NeiWiWo19,KhiSteiWi22_JCP,Hehl2022,KhiSteiWi23_JOTA}.
In our case the forward problem constitutes a second
(lower level) Nonlinear Optimization Problem (NLP),
i.e., a phase-field fracture NLP.
The objective of that lower level NLP consists in minimizing the
energy of the crack as it was formulated in the pioneering work \cite{FraMar98}.

The theoretical derivations on optimization in Banach spaces that we employ are based on
\cite{Maurer_Zowe:1979,Maurer:1981,Zowe_Kurcyusz:1979}. Applying these methods to a continuous
space-time phase-field fracture model and rigorously deriving the KKT conditions is novel.
Further we introduce these KKT conditions as a lower level problem within
the upper level optimal control problem and prove its regularity under certain conditions.

Concerning the mathematical analysis,
well-posedness with existence and convergence of
quasi-static brittle fracture
settings was investigated in \cite{FraLa03},
and in nonlinear elasticity in \cite{MaFraTo05}. In the year $2020$, fracture nucleation
was revisited as this is heavily used in numerous applications, but less clear from
a theoretical viewpoint \cite{KUMAR2020104027}.
Related works governing the analysis of phase-field fracture and damage models are
\cite{KneRoZa13, Lazzaroni2018, MIELKE20101242,PaThoTorWeiWie23,ThoBiWei20}.
Phase-field fracture for finite stresses was analysed in \cite{ThoBiWei20}. In particular,
the convergence of time-discrete solutions to solutions of the time-continuous problem was investigated.
The relationships between gradient damage, rate-independent damage and phase-field fracture are
discussed in \cite{BoVe16,PhaAmMarMau11} and \cite{Lazzaroni2018}.
For rate-independent damage models, \cite{KneRoZa13} introduced a vanishing viscosity approach.
Additionally, \cite{MIELKE20101242} emphasized complete (quasi-static) damage in particular,
avoiding the use of displacement fields in favour of stress and energy terms.
A monograph and a book chapter on rate independent systems and damage models
are available in \cite{Mielke2015,Mielke05}, respectively.

In summary, the main novelties in this work are:
The rigorous derivation of phase-field fracture as a complementarity system,
including the careful selection of suitable function spaces that are required
for the regularity of the resulting problems. Design of a space-time formulation of forward
phase-field fracture posed as complementarity system involving KKT conditions in Banach spaces rather than
a penalized space time-formulation \cite{KhiSteiWi22_JCP,KhiSteiWi23_JOTA}.
Here, corresponding cones and regularity results are rigorously investigated,
followed by second order necessary conditions and sufficient conditions.
Then, the upper level optimal control NLP with tracking type cost functional
and with the phase-field fracture KKT system as a constraint is studied
and KKT conditions, regularity properties, and further optimality conditions
are rigorously established.

The outline of this paper is as follows: In Section 2, the function
spaces and the basic notation are
introduced. In addition, the required abstract NLP theory is presented.
The phase-field fracture forward model is then introduced, along with its
formulation as an abstract NLP, in Section 3. Regular points of
that lower level NLP and its KKT system are given specific emphasis.
In Section 4, second order optimality conditions for the phase-field NLP
are presented. Then in Section 5, the upper level optimal control problem is formulated, where
the constraints are given by the previously derived
phase-field fracture KKT system.
Further, regular points of the optimal control NLP are characterized,
and its KKT conditions are derived.
Our work is summarized in Section 6.

\section{Theory for optimization in Banach spaces}

In this section, we introduce our notation,
recapitulate known results in the area of functional analysis from the literature,
and define function spaces required for space-time phase-field fracture.

\subsection{Notation}

Let $\Omega \subset \R^2$ be a bounded Lipschitz domain with boundary
partitioned as $\pd \Omega = \Gamma_N \cup \Gamma_D \cup \Gamma_F$,
where both $\Gamma_D$ and $\Gamma_N$ have nonzero one dimensional Hausdorff measure
and $\Omega \cup \Gamma_N$ is regular in the sense of Gröger \cite{Groeger1989}.
Next, define Hilbert spaces
$Q \define L^2(\Gamma_N)$ for the control force $q$,
$\Vphi \define H^2(\Omega) \hookrightarrow \Linf(\Omega)$
for the phase-field $\varphi$,
$\Vu \define H^1_D(\Omega; \R^2)$
for the two dimensional displacement field $u$, where
\begin{math}
  H^1_D(\Omega; \R^2) \define
  \defset{v \in H^1(\Omega; \R^2)}{v|_{\Gamma_D} = 0},
\end{math}
and the product space $V \define \Vu \times \Vphi$.
Finally consider a compact time interval $I \define [0,T]$ and define the spaces
\begin{equation}
  Y \define
  \Yu \times \Yphi \define
  L^2(I,\Vu) \times H^1(I,\Vphi)
  \qtextq{ and }
  \CIQ \define L^2(I, Q).
\end{equation}
We denote natural scalar products and norms with their space as index,
such as $\iprod[\Vu]{}{} \equiv \iprod[H_D^1(\Omega)]{}{}$ or
$\norm[\Yphi]{} \equiv \norm[H^1(I,H^2(\Omega))]{}$,
and $L^2$ scalar products and norms with their domain,
such as $\prodO{}{} \define \iprod[L^2(\Omega)]{}{}$
or $\prodIO{}{} \define \iprod[L^2(I, L^2(\Omega))]{}{}$.
Consequently the norm on $V$ is given by
\[
  \norm{\boldu}_V^2
  = \iprod[V]{\boldu}{\boldu}
  = \iprod[\Vu]{u}{u} + \iprod[\Vphi]{\varphi}{\varphi}
  = \sum_{\card\alpha \le 1} \prodO{D^\alpha u}{D^\alpha u}
  + \sum_{\card\alpha \le 2} \prodO{D^\alpha\varphi}{D^\alpha\varphi}.
\]
The norm on $Y$ is given by
\begin{align*}
  \norm{\boldu}_{Y}^2 = \int_I \bigl(
  \norm[\Vu]{u(t)}^2 + \norm[\Vphi]{\varphi(t)}^2 + \norm[\Vphi]{\dphi(t)}^2
  \bigr) \dt.
\end{align*}
Occasionally we will need the norm of a component
$\varphi \in H^1(I, \Vphi)$ or $u \in L^2(I, \Vu)$
and not of the combined function $\boldu = (u, \varphi)$:
\begin{align*}
  \norm[\Yu]{u}^2 &= \int_I \norm[\Vu]{u(t)}^2 \!\dt, &
  \norm[\Yphi]{\varphi}^2 &= \int_I \bigl(
  \norm[\Vphi]{\varphi(t)}^2 + \norm[\Vphi]{\dphi(t)}^2 \bigr) \dt.
\end{align*}

\subsection{Abstract NLP theory}

We consider the following
constrained nonlinear optimization problem of NLP type,
which was studied by Maurer and Zowe
\cite{Maurer_Zowe:1979,Maurer:1981}
based on prior work by Robinson and Kurcyusz
\cite{
  Robinson:1976b,Kurcyusz:1976,Zowe_Kurcyusz:1979}:
given Banach spaces $Y, Z$,
a closed convex set $C \subseteq Y$, 
a closed convex cone $K \subset Z$,
a cost functional $f\colon Y \to \R$,
and a constraint map $g\colon Y \to Z$,
the problem reads
\begin{equation}
  \label{eq:NLP}
  \min_{y \in C} \ f(y) \qstq g(y) \in K.
\end{equation}
As always in nonconvex optimization,
we regard every local minimizer as a solution,
and we are interested in first and second order conditions
that characterize these local minimizers.

The feasible set of \eqref{eq:NLP} is $M \define C \cap \Inv g(K)$.
For any feasible point $\by \in M$,
the \emph{tangent cone} $T(M, \by)$ and
the \emph{linearized cone} $L(M, \by)$ are defined as
\begin{align*}
  T(M, \by)
  &\define \defset{h \in C_{\by}}
    {h = \lim (y_n - \by) / t_n, \, y_n \in M, \, t_n > 0, \, t_n \to 0}, \\
  L(M, \by)
  &\define \defset{h \in C_{\by}}{g'(\by) h \in K_{g(\by)}}
    = C_{\by} \cap \Inv{g'(\by)}(K_{g(\by)}).
\end{align*}
Here we assume that the Fréchet derivatives
$f'(\by) \in L(Y, \R) = Y^*$ and $g'(\by) \in L(Y, Z)$ at $\by$ exist,
and by $C_{\by} \define \cone(C - \by)$
and $K_{g(\by)} \define \cone(K - g(\by))$ we
denote the conical hulls of $C - \by$ and $K - g(\by)$, respectively,
with $\cone(S) \define \defset{\alpha s}{s \in S, \, \alpha \ge 0}$.
By $Y^*$ we denote the topological dual space of $Y$,
and for any nonempty subset $S \subseteq Y$ we consider the \emph{dual cone}
$S^* \define \defset{l \in Y^*}{l s \ge 0 \tfor s \in S}$.
\begin{definition}
  A feasible point $\by \in M$ is called
  \emph{regular for \eqref{eq:NLP} (in the sense of Zowe and Kurcyusz)}
  if $g'(\by) C_{\by} - K_{g(\by)} = Z.$
\end{definition}
An important consequence of regularity is the inclusion
$L(M, \by) \subseteq T(M, \by)$.
KKT type optimality conditions for local minimizers of \eqref{eq:NLP}
are now given in the following theorem \cite{Maurer_Zowe:1979,Maurer:1981}
in terms of the Lagrangian $\Lagr\colon Y \times Z^* \to \R$ with
$\Lagr(y, l) \define f(y) - l g(y)$.
\begin{theorem}\label{KKT}
  If $\by \in M$ is a minimizer,
  then $f'(\by) h \ge 0$ for all $h \in T(M, \by)$,
  i.e., $f'(\by) \in T(M, \by)^*$.
  If, additionally, $\by$ is regular,
  then $f'(\by) h \ge 0$ for all $h \in L(M, \by)$,
  i.e., $f'(\by) \in L(M, \by)^*$.
  Equivalently, there is a Lagrange multiplier $l \in K^*$ such that
  \begin{equation}
    \label{eq:KKT}
    \pd_y \Lagr(\by, l) \equiv f'(\by) - lg'(\by) \in C_{\by}^*
    \qtextq{and}
    l g(\by) = 0.
  \end{equation}
\end{theorem}
In order to solve \eqref{eq:NLP} we have to find a \emph{KKT point}, i.e.,
a regular solution $\by$ of the KKT conditions \cref{eq:KKT}.
Assuming that second order Fréchet derivatives
$f''(\by) \in L(Y, L(Y, \R)) \cong L(Y, Y; \R)$ and
$g''(\by) \in L(Y, L(Y, Z)) \cong L(Y, Y; Z)$ exist
and that $C = Y$ (or $\by \in \tint C$),
hence $C_{\by} = Y$ and $C_{\by}^* = \set{0}$,
\cite{Maurer_Zowe:1979,Maurer:1981} also provide
second order necessary optimality conditions
at a given KKT point $(\by, l)$.
These conditions are stated in terms of the cones
$K^l \define K \cap \ker l$ and $T(M^l, \by)$
with $M^l \define \Inv g(K^l)$ as well as
\begin{equation*}
  L(M^l, \by) \define \defset{h \in Y}{g'(\by) h \in K^l_{g(\by)}}
  \qtextq{with}
  K^l_{g(\by)} \define K_{g(\by)} \cap \ker l.
\end{equation*}
The cone $L(M^l, \by)$ is called the \emph{critical cone} at $(\by, l)$.
\begin{theorem}
  \label{NSO}
  If $\by \in M$ is a minimizer
  with an associated Lagrange multiplier~$l$, then
  \begin{equation*}
    \pd_{yy} \Lagr(\by, l)(h, h) \ge 0 \qfor h \in T(M^l, \by).
  \end{equation*}
  If, additionally, $\by$ is \emph{regular w.r.t.\ $K^l$},
  that is, $g'(\by) Y - K^l_{g(\by)} = Z$, then
  \begin{equation*}
    \pd_{yy} \Lagr(\by, l)(h, h) \ge 0 \qfor h \in L(M^l, \by).
  \end{equation*}
\end{theorem}
In order to formulate sufficient optimality conditions,
the linearized cone $L(M, \by)$ has to be a good approximation
of the admissible set $M$ in a neighbourhood of $\by$.
\begin{definition}
  \label{def:approx-property}
  We say that the linearized cone $L(M, \by)$ \emph{approximates} $M$ at $\by$
  if there exists a map $h\colon M \to L(M, \by)$ such that
  \begin{equation*}
    \norm{h(y) - (y - \by)}  \in o(\norm{y - \by}) \qtextq{as} M \ni y \to \by.
  \end{equation*}
\end{definition}
Using the notation
$\bar B^Y_r(\by) \define \defset{y \in Y}{\norm{y - \by} \le r}$,
sufficient optimality conditions of first and second order
are now given in the following theorem \cite{Maurer_Zowe:1979,Maurer:1981}.
\begin{theorem}
\label{SFSO}
  Assume that $L(M, \by)$ approximates $M$ at $\by$.
  If there exists $\gamma > 0$ such that
  \begin{equation*}
    f'(\by) h \ge \gamma \norm{h} \qfor h \in L(M, \by),
  \end{equation*}
  then there exist $\alpha > 0$ and $\delta > 0$ with
  \begin{equation*}
    f(y) \ge f(\by) + \alpha \norm{y - \by}
   \qfor y \in M \cap \bar B^Y_\delta(\by).
  \end{equation*}
  If $(\by, l)$ is a KKT point
  and there exist $\gamma > 0$ and $\beta > 0$ such that
  \begin{equation*}
    \pd_{yy} \Lagr (\by, l)(h, h) \ge \gamma \norm{h}^2
    \qfor
    h \in L(M, \by) \cap \defset{h \in C_{\by}}{l g'(\by) h \le \beta \norm{h}},
  \end{equation*}
  then there exist $\alpha > 0$ and $\delta > 0$ with
  \begin{equation*}
    f(y) \ge f(\by) + \alpha \norm{y - \by}^2
    \qfor y \in M \cap \bar B^Y_\delta(\by).
  \end{equation*}
\end{theorem}

\section{Phase-field fracture as nonlinear energy minimization problem}
In this section, we formulate phase-field fracture as
an energy minimization problem with the crack irreversibility as
inequality constraint and the initial condition as equality constraint.
We formulate the corresponding cones and establish regularity
results and the KKT conditions.

\subsection{Problem statement}

Our phase-field fracture formulation differs from most works
found in the literature in that the crack irreversibility condition
is treated in a continuous-time fashion \cite{KhiSteiWi22_JCP, KhiSteiWi23_JOTA}
rather than by an (incremental) discrete-time formulation
\cite{BourFraMar00,Bour07,BuOrSue10,AmGeraLoren15,MieWelHof10a,KuMue10}.
We refer to the introduction for the discussion of various possibilities
of imposing the irreversibility constraint.
The continuous-time constraint enables us to formulate phase-field fracture
in a space-time setting, which we continue to utilize in this work.
The forward problem NLP$_E$ reads as follows:

\begin{Problem}\label{prob:NLP_E}
  Given the phase-field regularization $ε>0$, the bulk regularization $κ\in(0,1)$, 
  the Lam\'e parameters $μ>0$, $λ>-\frac23μ$, the critical energy release rate $G_c$,
  an initial phase-field $\varphi_0 \in \Vphi$
  and a space-time control $q \in \CIQ$,
  find a function $\boldu \in Y = \Yu \times \Yphi$
  consisting of a displacement field
  \begin{equation*}
    u = (u_1, u_2) \in \Yu = L^2(I, \Vu)
  \end{equation*}
  and a phase-field $\varphi \in \Yphi = H^1(I, \Vphi)$
  that minimize the crack energy $f\colon Y \to \R$ subject to
  the initial condition and the crack irreversibility constraint:
  \begin{equation}
    \label{eq:NLP_E}
    \begin{aligned}
      \min_{\boldu \in Y} \quad & f(\boldu) \define
      \frac12 \prodIO{\gkappa(\varphi) \C e(u)}{e(u)} + {} \\
      &\kern3.8em \frac{G_c}{2 \varepsilon} \normIO{1-\varphi}^2 +
      \frac{G_c \varepsilon}{2} \normIO{\nabla\varphi}^2 - \prodIN{q}{u} \\
      \stq
      &\kern0.6em \varphi(0) = \varphi_0 \qtextq{in} \Vphi, \\
      &{-}\dphi(t) \ge 0 \qtextq{in} \Vphi \aein I.
    \end{aligned}
  \end{equation}
\end{Problem}
In the integrand of the first term
we use the Frobenius scalar product $\C e(u) : e(u)$
of matrices $\C e(u)$ and $e(u)$ in $\R^{2 \times 2}$,
where $e(u) \define \frac{1}{2}(\nabla u + \nabla u^T)$
denotes the symmetric strain gradient
and $\C e(u) \define 2\mu e(u) + \lambda \text{tr}(e(u))I_2$ the stress tensor
with the identity matrix $I_2\in\R^{2 \times 2}$.
This product is multiplied with the nonlinear degradation function
$\gkappa(\varphi) = (1-\kappa) \varphi^2 + \kappa$.

\subsection{Constraints and cones}
In order to formulate the KKT system that corresponds to \eqref{eq:NLP_E},
we have to define
the operator $g$ (which is affine linear in this case),
the spaces $Z$, $Z^*$ and the cones $K$, $K^*$.
There is no set constraint, i.e., $C = Y$.
First we collect the equality and inequality constraints
$g_{\setE}$ and $g_{\setI}$, respectively,
to define $g\colon Y \to Z$ as
\begin{equation*}
  g(\boldu) =
  \begin{pmatrix}
    g_{\setE}(\boldu) \\
    g_{\setI}(\boldu)
  \end{pmatrix}
  \define
  \begin{pmatrix}
    \varphi(0) - \varphi_0 \\
    -\dphi
  \end{pmatrix}
  \qfor \boldu\in Y
  .
\end{equation*}
It is clear that the upper term belongs to $\Vphi$
for each $\varphi_0 \in \Vphi$
since $\varphi(t) \in \Vphi$ for all $t \in I$,
particularly $\varphi \in H^1(I, \Vphi) \hookrightarrow  C(I, \Vphi)$.
Moreover, the point evaluation $φ\mapsto\varphi(0)$
is surjective onto $\Vphi$.
By definition, the second component always belongs to $L^2(I, \Vphi)$.
Consequently we define the image space $Z$ as
\begin{equation*}
  Z \define Z_1 \times Z_2 \qtextq{with}
  Z_1 \define \Vphi, \quad
  Z_2 \define  L^2(I, \Vphi).
\end{equation*}
The dual space is
\begin{equation*}
  Z^* = Z_1^* \times Z_2^* = \Vphi^* \times L^2(I, \Vphi^*).
\end{equation*}
Next we define the constraints cone $K \subset Z$ as
\begin{equation*}
  K \define K_1 \times K_2 \qtextq{with}
  K_1 \define \set{0} \subset Z_1, \quad
  K_2 \define \defset{z_2 \in Z_2}{z_2 \ge 0}.
\end{equation*}
More precisely, $K_2$ has to be understood as
\begin{equation*}
  K_2 = \Defset{z_2 \in Z_2}
  {z_2(t) \ge 0 \text{ in } \Vphi \text{ for a.e.\ } t \in I}.
\end{equation*}
It is clear that $K_2$ and hence $K$ are closed convex cones.
Finally we need the dual cone $K^* \subset Z^*$.
Since $K = \set{0} \times K_2$, $K^*$ has the product structure
\begin{align*}
  K^*
  &= Z_1^* \times K_2^*
    = Z_1^* \times \defset{l_2 \in Z_2^*}{l_2 z_2 \ge 0 \tfor z_2 \in K_2} \\
  &= Z_1^* \times \Defset{l_2 \in L^2(I, \Vphi^*)}
    {\Idualprod{l_2(t)}{z_2(t)} \ge 0 \tfor z_2 \in K_2},
\end{align*}
where $\dualprod{u}{v}$ denotes the dual pairing between $\Vphi^*$ and $\Vphi$,
i.e., between $H^2(\Omega)^*$ and $H^2(\Omega)$.
Given a solution candidate $\bbu \in M$ with feasible set
\begin{equation*}
  M = \defset{\boldu = (u, \varphi) \in Y}{g(\boldu) \in K} =
  V_u \times \defset{\varphi \in \Vphi}{\varphi(0) = \varphi_0, \, -\dphi \ge 0},
\end{equation*}
we have $\bphi(t) = \bphi_0 + \int_0^t \dot\bphi(s) \ds \le \bphi_0$ for all $t\in I$ and
the relevant cones are
\begin{align}
  \label{eq:CandK}
  C_{\bbu} &= Y,
  & C_{\bbu}^* &= \set{0},
  & K_{g(\bbu)}
  &= \set{0} \times (K_2 + \Span\set{\dot\bphi}).
\end{align}
Next we have to compute the derivatives
$f'(\bbu) \in Y^*$ and $g'(\bbu) \in L(Y, Z)$.
For $f'(\bbu)$ see the following proposition.
Since $g$ is affine linear,
for any given direction $\boldPhi \define (\Phi_u, \Phiphi) \in Y$ we obtain
\begin{equation*}
  g'(\bbu)(\boldPhi) =
  \begin{pmatrix}
    \Phiphi(0) \\
    -\dPhiphi
  \end{pmatrix}.
\end{equation*}
Finally, given any multiplier $\boldl = (l_1, l_2) \in K^*$,
we have $\boldl g'(\bbu) \in Y^*$ with
\begin{equation*}
  \boldl g'(\bbu)(\boldPhi)
  =
  \dualprod{l_1}{\Phiphi(0)} - \Idualprod{l_2(t)}{\dPhiphi(t)}.
\end{equation*}
In preparation of some estimates to be employed in the proof of \cref{prop:f'+f''},
we relate certain expressions to the norm in $Y$:
\begin{lemma}\label{lem:estimatescalarproduct}
Let $μ>0$, $λ>-\frac23 μ$. Then there is $C>0$ such that the following
estimates hold for all $\boldu=(u,φ)\in Y$, $\boldu_1,\boldu_2\in Y$:
\begin{align*}
  \prodIO{ψ \C e(u_1)}{e(u_2)}
  &\le C \norm[\Linf(\IO)]{ψ} \norm[Y]{\boldu_1} \norm[Y]{\boldu_2}
    \qfor ψ \in \Linf(\IO), \\
 \norm[\Linf(\IO)]{φ} &\le C \norm[Y]{\boldu}, \\
 \max\set{\normIO{∇φ}, \normIO{φ}} &\le \norm[Y]{\boldu}.
\end{align*}
\end{lemma}
\begin{proof}
The last inequality is straightforward, as
\begin{align*}
  \normIO{∇φ}
  &\le \norm[L^2(I,H^2(\Omega))]{φ} \le \norm[Y]{\boldu}, \\
  \normIO{φ}
  &\le \norm[L^2(I,H^2(\Omega))]{φ} \le \norm[Y]{\boldu},
\end{align*}
whereas the second estimate holds due to the embedding
$\Yphi = H^1(I, H^2(\Omega)) \hookrightarrow \Linf(I, \Linf(\Omega))$.
The first one mainly relies on Hölder's inequality:
With some $C>0$ (depending on $μ$ and~$λ$),
\begin{align*}
  \prodIO{ψ \C e(u_1)}{e(u_2)}
  &= \int_{\IO} ψ(t,x) \C e(u_1(t,x)) : e(u_2(t,x)) \dx \dt \\
  &\le C \norm[\Linf(\IO)]{ψ} \norm[\IO]{∇u_1} \norm[\IO]{∇u_2} \\
  &\le C \norm[\Linf(\IO)]{ψ} \norm[Y]{\boldu_1} \norm[Y]{\boldu_2}.
\end{align*}
\end{proof}
\newcommand\nn{\nonumber}
\begin{proposition}
\label{prop:f'+f''}
The energy functional $f$ defined in \eqref{eq:NLP_E}
is twice Fréchet differentiable in $Y$.
For any direction $\boldPhi \in Y$ the first and second derivatives
at any element $\boldu \in Y$ are given by
\begin{equation}
\label{eq:derivatives}
\begin{aligned}
  f'(\boldu)(\boldPhi)
  &= \prodIO{\gkappa(\varphi) \C e(u)}{e(\Phiu)} - \prodIN{q}{\Phiu}
    + G_c \varepsilon \prodIO{\nabla\varphi}{\nabla\Phiphi} \\
  &\quad- \frac{G_c}{\varepsilon} \prodIO{1 - \varphi}{\Phiphi}
    + (1 - \kappa) \prodIO{\varphi \Phiphi \C e(u)}{e(u)}, \\
  f''(\boldu)(\boldPhi,\boldPhi)
  &= \prodIO{\gkappa(\varphi)\C e(\Phiu)}{e(\Phiu)}
    + 4 (1 - \kappa) \prodIO{\varphi\Phiphi \C e(u)}{e(\Phiu)} \\
  &\quad+ G_c \varepsilon \normIO{\nabla\Phiphi}^2
    + \frac{G_c}{\varepsilon} \normIO{\Phiphi}^2
    + (1 - \kappa) \prodIO{\Phiphi^2 \C e(u)}{e(u)}.
\end{aligned}
\end{equation}
\end{proposition}

\newcommand\boldPsi{\bs \Psi}
\newcommand\Psiu[1][]{\Psi_{u\optindex{#1}}}
\newcommand\Psiphi[1][]{\Psi_{\varphi\optindex{#1}}}

\begin{proof}
The expressions in \eqref{eq:derivatives} are easily confirmed to be Gâteaux derivatives.
Moreover, for any $\boldu,\boldPhi,\boldPsi\in Y$, by
\cref{lem:estimatescalarproduct} in combination with the identity
$\gkappa(\varphi+\Psiphi)=\gkappa(\varphi)+(1-κ)(2φ\Psiphi+\Psiphi^2)$
we have
\begin{align*}
  \abs{f'
  &(\boldu+ \boldPsi) \boldPhi - f'(\boldu)\boldPhi} \\
  &= \Bigl|\prodIO{\gkappa(φ+\Psiphi)\C e(\Psiu)}{e(\Phiu)} \\
  &\quad\;+ (1-κ)\prodIO{(2φ\Psiphi+\Psiphi^2)\C e(u)}{e(\Phiu)}
    + G_cε \prodIO{∇\Psiphi}{∇\Phiphi}\\
  &\quad\;+ \frac{G_c}{ε}\prodIO{\Psiphi}{\Phiphi}
    + (1-κ) \prodIO{\Psiphi\Phiphi\C e(u+\Psiu)}{e(u+\Psiu)}\\
  &\quad\;+ 2(1-κ)\prodIO{φ\Phiphi\C e(u)}{e(\Psiu)}
    + (1-κ)\prodIO{φ\Phiphi\C e(\Psiu)}{e(\Psiu)}\Bigr|\\
  &\le C\norm[\Linf(\IO)]{\gkappa(φ+\Psiphi)} \norm[Y]{\boldPsi} \norm[Y]{\boldPhi} \\
  &\quad+ (1-κ) C\norm[\Linf(\IO)]{2φ\Psiphi+\Psiphi^2}
    \norm[Y]{\boldu} \norm[Y]{\boldPhi}
    + G_cε \norm[Y]{\boldPsi}\norm[Y]{\boldPhi} \\
  &\quad+ \frac{G_c}{ε}\norm[Y]{\boldPsi}\norm[Y]{\boldPhi}
    + (1-κ)C \norm[\Linf(\IO)]{\Psiphi\Phiphi}\norm[Y]{\boldu+\boldPsi}^2\\
  &\quad+ 2(1-κ)C\norm[\Linf(\IO)]{φ\Phiphi}
    \norm[Y]{\boldu}\norm[Y]{\boldPsi}
    + (1-κ)C\norm[\Linf(\IO)]{φ\Phiphi}\norm[Y]{\boldPsi}^2,
\end{align*}
where $C>0$ is the constant from \cref{lem:estimatescalarproduct}, another application of which yields
\begin{align*}
  \norm[Y^*]{f'(\boldu+\boldPsi) - f'(\boldu)}
  &\le C\norm[\Linf(\IO)]{\gkappa(φ+\Psiphi)} \norm[Y]{\boldPsi}  \\
  &\quad+ (1-κ) C\norm[\Linf(\IO)]{2φ\Psiphi+\Psiphi^2}
    \norm[Y]{\boldu} + G_cε \norm[Y]{\boldPsi} \\
  &\quad+ \frac{G_c}{ε}\norm[Y]{\boldPsi}
    + (1-κ)C^3 \norm[Y]{\boldPsi}\norm[Y]{\boldu+\boldPsi}^2\\
  &\quad+ 2(1-κ)C^3\norm[Y]{\boldu}^2\norm[Y]{\boldPsi}
    + (1-κ)C^3\norm[Y]{\boldu}\norm[Y]{\boldPsi}^2.
\end{align*}
Hence, apparently, $f'(\boldu+\boldPsi)\to f'(\boldu)$ in $Y^*$ as $\boldPsi\to 0$.
Continuity of the Gâteaux derivative implies Fréchet differentiability (cf. \cite[Prop.~4.8(c)]{zeidlerI}).

For the second derivative, we similarly obtain for any $\boldu,\boldPhi,\boldPsi\in Y$
\begin{align*}
  \abs{f''
  &(\boldu+\boldPsi)(\boldPhi,\boldPhi)-f''(\boldu)(\boldPhi,\boldPhi)} \\
  &= \Bigl|(1-κ)\prodIO{(2φ\Psiphi+\Psiphi^2)\C e(\Phiu)}{e(\Phiu)}
    + 4 (1 - \kappa) \prodIO{\varphi\Phiphi \C e(\Psiu)}{e(\Phiu)}\nn \\
  &\quad\;+ 4 (1 - \kappa) \prodIO{\Psiphi\Phiphi \C e(u+\Psiu)}{e(\Phiu)}\nn
    + 2(1 - \kappa) \prodIO{\Phiphi^2 \C e(u)}{e(\Psiu)} \\
  &\quad\;+ (1 - \kappa) \prodIO{\Phiphi^2 \C e(\Psiu)}{e(\Psiu)}\Bigr| \\
  &\le (1-κ)C \norm[\Linf(\IO)]{2φ\Psiphi+\Psiphi^2}\norm[Y]{\boldPhi}^2
    + 4(1-κ)C^3\norm[Y]{\boldu}\norm[Y]{\boldPsi}\norm[Y]{\boldPhi}^2\\
  &\quad+ 4(1 - \kappa)C^3\norm[Y]{\boldPsi}
    \norm[Y]{\boldu+\boldPsi}\norm[Y]{\boldPhi}^2
    + 2(1 - \kappa)C^3\norm[Y]{\boldPhi}^2\norm[Y]{\boldu}\norm[Y]{\boldPsi}\\
  &\quad+ (1-κ)C^3\norm[Y]{\boldPhi}^2\norm[Y]{\boldPsi}^2,
\end{align*}
so that
\begin{align*}
  \norm[L(Y,Y;ℝ)]{f''(\boldu+\boldPsi)-f''(\boldu)}
  &\le (1-κ)C \norm[\Linf(\IO)]{2φ\Psiphi+\Psiphi^2} \\
  &\quad+ 6(1-κ)C^3\norm[Y]{\boldu}\norm[Y]{\boldPsi} \\
  &\quad+ 4(1 - \kappa)C^3\norm[Y]{\boldPsi}\norm[Y]{\boldu+\boldPsi}
    + (1-κ)C^3\norm[Y]{\boldPsi}^2.
\end{align*}
Once more continuity of the Gâteaux derivative ensures Fréchet differentiability.
\end{proof}

\subsection{Regularity and KKT system}
\label{sec:regularity}
Recall that a feasible point $\bbu \in M$
is regular for \eqref{eq:NLP_E} if
\[ 
  Z = g'(\bbu) Y - K_{g(\bbu)},
\] 
where in our case
\begin{align*}
  Z &= Z_1 \times Z_2 = \Vphi \times L^2(I, \Vphi), \\
  g'(\bbu) Y &= \Defset{\begin{pmatrix}
      \Phiphi(0) \\
      -\dPhiphi
    \end{pmatrix}}{(\Phiu, \Phiphi) \in Y}, \\
  K_{g(\bbu)} &= \Defset{\begin{pmatrix}
      0 \\
      k_2
    \end{pmatrix}
  + \alpha
  \begin{pmatrix}
    0 \\
    \dot\bphi
  \end{pmatrix}}{k_2 \in K_2 \text{ and } \alpha \in \R}
  .
\end{align*}
In fact, \emph{every} feasible point is regular
since $Z = g'(\bbu) Y$ by the following result.
\begin{lemma}
  \label{prop:surj}
  The derivative $g'(\bbu)\colon Y \to Z$ is surjective for every $\bu \in M$.
\end{lemma}
\begin{proof}
  Given $\bbu \in M$ and $z = (z_1, z_2) \in Z$,
  simply choose any $\Phiu \in \Yu$ and set
  $\Phiphi(t) \define z_1 - \int_0^t z_2(s) \ds$
  to obtain $\boldPhi \in Y$ with $g'(\bbu) \boldPhi = z$.
\end{proof}

  We work with $\varphi \in H^1(I,H^2(\Omega))$ to make use of the embedding
  $H^2(\Omega) \hookrightarrow \Linf(\Omega)$
  in proving \cref{prop:f'+f''}.
  A seemingly natural alternative would be to require less spatial regularity,
  for example by merely assuming $φ(t)\in H^1(\Omega)$ for every $t\in I$. In this case, however,
  it would be unclear whether $f$ was well-defined as a functional $f\colon Y\to \R$,
  because finiteness of the term $\prodO{\gkappa(\varphi) \C e(u)}{e(u)}$ could not be guaranteed
  for each $\boldu \in Y$. Of course, it would be possible to replace $\gkappa$ by a bounded (cut off) variant,
  but even then, related terms (and thus the same problem) would re-emerge when dealing with derivatives.
  Another common space for $\varphi$ is $L^2(I, H^1(\Omega))$ with
  $\dphi \in L^2(I, H^1(\Omega)^*)$.
  The following counterexample shows that this choice of function spaces
  cannot ensure that any feasible point is regular for \eqref{eq:NLP_E}.
  \begin{lemma}
  Let $I=[0,1]$ and $\Omega=(0,1)^2$, let $\bphi$ be a nonnegative function
  in $L^2(I,H^1(\Omega))$ such that
  $\dot \bphi \in L^2(I,H^1(\Omega)^*)$ with $\dot \bphi \le 0$.
  Then there exists a function $z\in V_{φ}\times L^2(I, H^1(\Omega)^*)$
  for which the equations $z_1=\Phiphi(0)$, $z_2 = -\dPhiphi + k_2 + \alpha \dot \bphi$
  cannot be fulfilled with any $\boldPhi\in L^2(I,H^1(\Omega))$
  that satisfies $\dPhiphi \in L^2(I,H^1(\Omega)^*)$,
  any $k_2 \in L^2(I, H^1(\Omega)^*)$ that satisfies $k_2\ge 0$, and any $\alpha \in \R$.
  \end{lemma}
  \begin{proof}
  Let $z_1 = \bphi(0) \in L^2(\Omega)$,
  and $z_2(t,x,y) = -x^{-\frac14}$ for every $t\in I$ and $(x,y)\in \Omega$.
  Then, for each $t\in I$, $z_2(t) \in L^2(\Omega) \subset H^1(\Omega)^*$ since
  \begin{equation*}
    \int_0^1 \int_0^1 \abs{z_2(t,x,y)}^2 \dx \diff y =
    \int_0^1 \int_0^1 \abs{-x^{-\frac14}}^2 \dx \diff y =
    \int_0^1 x^{-\frac12} \dx = 2 < \infty.
  \end{equation*}
  Moreover, this shows $z_2\in L^2(I,H^1(\Omega)^*)$.
  Now we seek functions $\Phiphi\in L^2(I,H^1(\Omega))$ satisfying
  $\dPhiphi \in L^2(I, H^1(\Omega)^*)$ and
  $k_2 \in L^2(I,H^1(\Omega)^*)$ with $k_2 \ge 0$ and constant $\alpha \in \R$ such that
  \begin{equation*}
    \Phiphi(0) = z_1\qtextq{and}
    \dPhiphi = -z_2 + k_2 + \alpha \dot \bphi.
  \end{equation*}
  If such a function $\Phiphi$
  exists, then for almost every $t \in I$ we have $\Phiphi(t)\in H^1(\Omega)$ and
  \begin{equation}\label{eq:counterexample-estimate}
    \Phiphi(t) =
    \int_0^t [k_2(s) - z_2(s) + \alpha \dot \bphi(s)] \ds
  \end{equation}
  in $H^1(\Omega)$. If $\alpha \le 0$, then due to $k_2 \ge 0$ and $\dot \bphi \le 0$ this ensures
  \[
   \Phiphi(t) \ge - \int_0^t z_2(s) \ds = -tz_2(0).
  \]
  If, on the other hand $\alpha > 0$, \cref{eq:counterexample-estimate} shows that
  $\Phiphi(t) \ge -t z_2(0) + \alpha \bphi(t) \ge -t z_2(0)$.
  In both cases this implies the existence
  of a function $g \in H^1(\Omega)$ with
  $g(x,y) \ge -z_2(0,x,y) = {} x^{-\frac14}$,
  namely $g(x,y) = \frac{1}{t}\Phiphi(t,x,y)$
  for some $t>0$ for which \eqref{eq:counterexample-estimate} is valid.
  Now let $\tilde{g}(x) = \int_0^1 g(x,y) \diff y$.
  Then $\tilde{g} \in H^1((0,1)) \subset C([0,1])$
  and $\tilde{g}(x) \ge x^{-\frac14}$ for all $x\in(0,1)$,
  which contradicts $\tilde{g} \in C([0,1])$
  since continuous functions are bounded.
  Therefore such functions $\Phiphi$ and $k_2$ cannot exist.
\end{proof}
Now we are ready to formulate the KKT system
that corresponds to \eqref{eq:NLP_E}.
\begin{proposition}
  \label{KKTE}
  Let the data and parameters be as in \cref{prob:NLP_E}.
  Since every local minimizer $\bbu = (\bu, \bphi) \in M$ of NLP$_E$
  is regular by \cref{prop:surj},
  there exists a multiplier $\boldl = (l_1, l_2) \in Z^*$
  such that the KKT conditions of \cref{KKT} hold:
  \begin{align}
    \label{KKT_1}\tag{KKT\,1}
    -\dot\bphi(t) &\ge 0 \text{ in } \Vphi \text{ for a.e.\ }t\in I,\\
    \label{KKT_2}\tag{KKT\,2}
    \bphi(0) &= \varphi_0 \text{ in } \Vphi, \\
    \label{KKT_3}\tag{KKT\,3}
    l_2(t) &\ge 0 \text{ in } \Vphi^* \text{ for a.e.\ }t\in I,
  \end{align}
  \vspace*{-1.9pc}
  \begin{align}
    \notag
    \prodIO{\gkappa(\bphi) \C e(\bu)}{e(P_u(\fcdot))} - \prodIN{q}{P_u(\fcdot)}
    + G_c \varepsilon \prodIO{\nabla\bphi}{\nabla P_\varphi(\fcdot)} \\
    \notag
    {} - \frac{G_c}{\varepsilon} \prodIO{1 - \bphi}{P_\varphi(\fcdot)}
    + (1 - \kappa) \prodIO{\bphi P_\varphi(\fcdot) \C e(\bu) }{e(\bu)} \\
    \label{KKT_4}\tag{KKT\,4}
    {} - \dualprod{l_1}{P_\varphi(\fcdot)(0)}
    + \Idualprod{l_2(t)}{\pd_t P_\varphi(\fcdot)(t)} &= 0, \\
    \label{KKT_5}\tag{KKT\,5}
    \dualprod{l_1}{\bphi(0) - \varphi_0} - \Idualprod{l_2(t)}{\dot\bphi(t)} &= 0,
  \end{align}
  where we introduce canonical projections $P_u(\boldPhi) = \Phiu$
  and $P_\varphi(\boldPhi) = \Phiphi$.
\end{proposition}
\begin{proof}
  Conditions \eqref{KKT_1} and \eqref{KKT_2}
  are just feasibility $\bbu \in M$.
  Condition \eqref{KKT_3} is equivalent to $\boldl \in K^*$.
  The stationarity condition $f'(\bbu) - \boldl g'(\bbu) \in C_{\bbu}^*$
  in the form \eqref{KKT_4} is immediate from the representation of $f'$
  in \cref{prop:f'+f''} since $C_{\bbu}^* = \set{0}$.
  Finally, \eqref{KKT_5} is the complementarity condition
  $\boldl g(\bbu) = 0$.
\end{proof}

\section{Further optimality conditions}
\label{sec_further_opt_cond}

Here we formulate the specific second order necessary conditions
as well as first and second order sufficient conditions
of the forward problem NLP$_E$.
\subsection{Necessary optimality conditions of second order}

The Lagrangian corresponding to \eqref{eq:NLP_E} reads
\begin{align*}
  \Lagr(\boldu,\boldl)
  &= f(\boldu) - \boldl g(\boldu) \\
  &= \frac12 \prodIO{\gkappa(\varphi) \C e(u)}{e(u)}
    + \frac{G_c}{2 \varepsilon} \normIO{1-\varphi}^2
    + \frac{G_c \varepsilon}{2} \normIO{\nabla\varphi}^2 - \prodIN{q}{u} \\
  &\quad- \dualprod{l_1}{\varphi(0) - \varphi_0} + \Idualprod{l_2(t)}{\dphi(t)}.
\end{align*}
By affine linearity of $g$ and \cref{prop:f'+f''},
the derivative $\pd_{\boldu \boldu}\Lagr(\bbu,\boldl)(\boldPhi,\boldPhi)$
for a pair $(\bbu, \boldl) \in Y \times Z^*$ and a
direction $\boldPhi = (\Phiu,\Phiphi) \in Y$ becomes
\begin{align*}
  \pd_{\boldu \boldu}\Lagr(\bbu,\boldl)(\boldPhi,\boldPhi)
  &= f''(\bbu)(\boldPhi,\boldPhi) \\
  &= \prodIO{\gkappa(\bphi)\C e(\Phiu)}{e(\Phiu)}
    + 4 (1 - \kappa) \prodIO{\bphi \Phiphi \C e(\bu)}{e(\Phiu)} \\
  &\quad+ G_c \varepsilon \normIO{\nabla\Phiphi}^2
    + \frac{G_c}{\varepsilon} \normIO{\Phiphi}^2
    + (1 - \kappa) \prodIO{\Phiphi^2 \C e(\bu)}{e(\bu)}.
\end{align*}
In order to ensure the required regularity of some minimizer $\bbu \in M$
with a multiplier $\boldl \in K^*$ for the second order necessary
optimality conditions, we have to show that
\begin{equation*}
  g'(\bbu) Y - K^{\boldl}_{g(\bbu)} = Z.
\end{equation*}
This is indeed the case for \emph{every} feasible point
since $g'(\bbu)\colon Y \to Z$ is surjective,
as shown in \cref{prop:surj}.
Finally we can apply \cref{NSO} to \cref{prob:NLP_E}:
\begin{proposition}
  Every local minimizer $\bbu$ of NLP$_E$
  is regular with respect to~$K^{\boldl}$,
  and for every
  \begin{math}
    \boldPhi \in L(M^{\boldl}, \bbu) =
    \defset{h \in Y}{g'(\bbu) h \in K^{\boldl}_{g(\bbu)}}
  \end{math}
  it holds
  \begin{align*}
    \prodIO{\gkappa(\bphi)\C e(\Phiu)}{e(\Phiu)}
    + 4 (1 - \kappa) \prodIO{\bphi \Phiphi \C e(\bu)}{e(\Phiu)} \\ {}
    + G_c \varepsilon \normIO{\nabla\Phiphi}^2
    + \frac{G_c}{\varepsilon}\normIO{\Phiphi}^2
    + (1 - \kappa) \prodIO{\Phiphi^2\C e(\bu)}{e(\bu)}
    &\ge 0.
  \end{align*}
\end{proposition}

\subsection{Sufficient optimality conditions}
First of all we note that the approximation property of \cref{def:approx-property}
is trivially satisfied at every feasible point $\bbu \in M$
since $M \subseteq L(M, \bbu)$. The sufficient optimality conditions \cref{SFSO}
of first and second order for \cref{eq:NLP_E} take the following form.
\begin{proposition}
\label{prop:1_2_sufficient}
Let $(\bbu, \boldl) \in Y \times Z^*$ be a KKT point for \cref{prob:NLP_E}.
The first order sufficient condition holds if
there exists $\gamma > 0$ such that
  \[
    f'(\bbu)(\boldPhi) \ge \gamma \norm[Y]{\boldPhi}
    \qfor \boldPhi \in L(M, \bbu),
  \]
  where
  \[
    L(M,\bbu) \define \defset{\boldPhi \in Y}
    {\Phiphi(0) = 0, \, -\dPhiphi \in K_2 + \Span\set{\dot\bphi}}.
  \]
  The second order sufficient condition holds if
  there exist $\gamma > 0$ and $\beta > 0$ such that
  \[
    \partial_{\boldu \boldu} \Lagr(\bbu,\boldl)(\boldPhi,\boldPhi) \ge
    \gamma \norm[Y]{\boldPhi}^2 \qfor \boldPhi \in L_\beta(M, \bbu),
  \]
  where
  \begin{align*}
    L_\beta(M, \bbu) \define
    \Bigl\{ \boldPhi \in Y\colon &\Phiphi(0) = 0, \,
             -\dPhiphi \in K_2 + \Span\set{\dot\bphi}, \Bigr. \\
    \Bigl. &{-}\int_I \dualprod{l_2(t)}{\dPhiphi(t)} \dt
      \le \beta \norm[Y]{\boldPhi} \Bigr\}
  \end{align*}
\end{proposition}
In \cref{lem:suff1} and \cref{lem:suff2} we show that
the sufficient optimality conditions of first and second order
do \emph{not} hold at any KKT point.
\begin{lemma}
 \label{lem:suff1}
 Let $(\bbu, \boldl) \in Y \times Z^*$ be a KKT point for \cref{prob:NLP_E}
 with $\bphi\not\equiv \varphi_0$.
 For each $\gamma >0$ there exists $\boldPhi \in L(M, \bbu)$ such that
 $f'(\bbu)(\boldPhi) {} = 0 < \gamma \norm[Y]{\boldPhi}$.
\end{lemma}
\begin{proof}
  Set $\Phiu = 0$ and $\Phiphi = -\varphi_0 + \bphi$.
  Then $\boldPhi \in L(M, \bbu)$ since $\Phiphi(0) = 0$ by \eqref{KKT_2}
  and $-\dPhiphi = -\dot \bphi \in K_2 + \Span\{\dot\bphi\}$.
  Further, due to \eqref{eq:derivatives},
  \eqref{KKT_4} and \eqref{KKT_5}, we get
  \begin{align*}
    f'(\bbu)(\boldPhi)
    &{}= \boldl g'(\bbu)
      = \dualprod{l_1}{\Phiphi(0)} - \Idualprod{l_2(t)}{\dPhiphi(t)} \\
    &=\dualprod{l_1}{\bphi(0)-\varphi_0} - \Idualprod{l_2(t)}{\dot\bphi(t)}
      = 0.
  \end{align*}
\end{proof}
\newcommand{\HC}{\mathcal{H}(\Omega)}

The following lemma serves as technical preliminary for the treatment of 
the sufficient second order condition.
It will be applied to $g=-\dot\bphi$ and its use could be entirely avoided under 
the additional assumption that $\dot\bphi$ be continuous.
\begin{lemma}
  \label{lem:Lusin_sandwich}
  Let $g\in L^2(I,\HC) {} \without{0}$ satisfy $g\ge 0$,
  where $\HC$ is a separable Banach space of functions
  with continuous embedding $\HC \hookrightarrow C^0(\Ombar)$.
  Then there are a subset $J\subset I$ of positive (Lebesgue) measure $\mu(J)$
  and a function $\Psi\in C_0^\infty(\Omega)\without{0}$
  such that $0\le \Psi(x)\le g(t,x)$ for almost every $(t,x)\in J\times\Omega$.
\end{lemma}
\begin{proof}
 As $g\not\equiv 0$, there is $ε>0$ such that the set $S_{ε} \define \set{(t,x)\in I\times\Omega \mid g(t,x)>2ε}$ 
 fulfills $\mu(S_{ε})>0$.
 We use Lusin's theorem (for a variant in Bochner spaces we refer to \cite{loeb_talvila}) to find a %
 set $K\subset I$ with $\mu(I\setminus K)<\mu(S_{ε}) / (2 μ(\Omega))$ such that $g\vert_{K}\colon K\to \HC$ 
 is continuous.
 Noting that the set $S\define (K\times\Omega)\cap S_{ε}$ has positive measure, we pick $(t_0,x_0)\in S$ 
 such that for all $δ>0$ we
 have $\mu(K\cap (t_0-δ,t_0+δ))>0$. 
  Relying on the continuity of $g(t_0,\fcdot)$, we pick $δ_1>0$ such 
  that for every $x\in B\define B_{δ_1}(x_0)$ we have $g(t_0,x)>2ε$,
  and from continuity of $g\vert_{K}$ and the embedding $\HC\hookrightarrow C^0(\Ombar)$
  we obtain some $δ_2>0$ such that for every $τ\in J\define K\cap (t_0-δ_2,t_0+δ_2)$
  we have $\norm[C^0(\Ombar)]{g(τ,\fcdot)-g(t_0,\fcdot)}<ε$. For every $τ\in J$ and $x\in B$ we thus have
 \[
  g(τ,x)=(g(τ,x)-g(t_0,x))+g(t_0,x) > -ε +2ε=ε.
 \]
 We finally choose any nonnegative $\Psi\in C^\infty_0(\Omega) \setminus\set{0}$ 
 which vanishes on $\Omega\setminus B$ and
 satisfies $\Psi(x)\in [0,ε]$ for every $x\in B$.
\end{proof}

\begin{lemma}
\label{lem:suff2}
 Let $(\bbu, \boldl) \in Y \times Z^*$ be a KKT point for \cref{prob:NLP_E}.
 Then for all $\gamma>0$ and $\beta >0$
 there exists $\boldPhi \in L_\beta(M, \bbu)$ such that
 $\partial_{\boldu \boldu} \Lagr(\bbu,\boldl)(\boldPhi,\boldPhi) <
    \gamma \norm[Y]{\boldPhi}^2$.
\end{lemma}
\begin{proof}
  Applying \cref{lem:Lusin_sandwich} to $-\dot \bphi$ with
  $\HC = \Vphi = {} H^2(\Omega) \hookrightarrow C^0(\Ombar)$
  leads to a set $J \subset I$ and $\Psi \in \Vphi$ such that
  $0 \le \Psi(x) \le -\dot \bphi(t,x)$
  for a.e.\ $(t,x)\in J\times\Omega$.
  For $\eta \in(0,μ(J))$ choose $J_\eta \subset J$
  with $\mu(J_\eta) = \eta$ and set
  $f_\eta(t) = \frac1{\eta} \int_0^t \chi_{J_\eta}(s) \ds$,
  where  $\chi_{J_\eta}$ denotes the characteristic function of $J_\eta$.
  We define $\boldPhi_{η}=(\Phiu,\Phiphi)$ by setting $\Phiu \equiv 0$
  and $\Phi_{\varphi}(t,x) = - f_\eta(t) \Psi(x)$. Then
  $\Phi_{\varphi}(0,x) = 0 \Psi(x) = 0$,
  $-\dot \Phi_{\varphi,\eta}(t,x) = \frac{1}{\eta}\chi_{J_\eta}(t) \Psi(x) \ge 0$,
  and \eqref{KKT_3} and \eqref{KKT_5} in conjunction with \eqref{KKT_2}
  yields
  \begin{align*}
    -\Idualprod{l_2(t)}{\dPhiphi(t)}
    &= \Idualprod{l_2(t)}{\tfrac{1}{\eta}\chi_{J_\eta}(t)\Psi}
      = \frac{1}{\eta} \int_{J_\eta}\dualprod{l_2(t)}{\Psi} \dt\\
    &\le \frac{1}{\eta}\int_{J_\eta} \dualprod{l_2(t)}{-\dot \bphi(t)}
      = 0\le \beta \norm[Y]{\Phi_{η}}.
\end{align*}
Therefore $\boldPhi_{η} \in  L_\beta(M, \bbu)$ for every $\beta > 0$. From
\begin{align*}
 \norm[Y]{\boldPhi_{η}}^2 &= \int_I (\norm[\Vphi]{\Phiphi(t)}^2 + \norm[\Vphi]{\dPhiphi(t)}^2) \dt
\ge \int_I \norm[\Vphi]{\dPhiphi(t)}^2 \dt \\
 &= \int_I \dot f_{η}(t)^2 \norm[\Vphi]{\Psi}^2 \dt
 = \frac{\norm[\Vphi]{\Psi}^2}{η^2} \int_I \chi_{J_{η}}(t)^2 \dt
 = \frac1{η}\norm[\Vphi]{\Psi}^2
\end{align*}
we can conclude $ \norm[Y]{\boldPhi_{η}}^2 \to \infty$ as $\eta \to 0$. With $c_1>0$
being the constant introduced in \cref{lem:estimatescalarproduct}, we obtain
\begin{align*}
  \partial_{\boldu \boldu} \Lagr(\bbu,\boldl)(\boldPhi_{η},\boldPhi_{η})
  &= G_c \varepsilon \normIO{\nabla\Phiphi}^2
    + \frac{G_c}{\varepsilon}\normIO{\Phiphi}^2
    + (1 - \kappa) \prodIO{\Phiphi^2 \C e(\bu)}{e(\bu)} \\
  &\le G_c \varepsilon \normIO{\nabla\Phiphi}^2
    + \frac{G_c}{\varepsilon}\normIO{\Phiphi}^2
    + c_1\norm[\Linf(I\times \Omega)]{\Phiphi}^2 \normIO{\bu}^2.
\end{align*}
As $\Phiphi(t)=-f_{η}(t)\Psi$ and $|f_\eta(t)|\le 1$ for all $t\in I$,
this estimate together with the embedding $\Vphi = {} H^2(\Omega) \hookrightarrow \Linf(\Omega)$ shows that
there exists a constant $C>0$  such that
$ \partial_{\boldu \boldu} \Lagr(\bbu,\boldl)(\boldPhi_{η},\boldPhi_{η}) {}\le C \norm[\Vphi]{\Psi}^2 $
for all $\eta{}\in(0,μ(J))$.
Consequently for every $\gamma >0$ we can find $\eta$ small enough such that
$ \partial_{\boldu \boldu} \Lagr(\bbu,\boldl)(\boldPhi_{η},\boldPhi_{η}) < \gamma \norm[Y]{\boldPhi_{η}}^2$.
\end{proof}
\section{Upper level NLP with phase-field constraint}
In this section we consider an optimal control NLP
whose constraints are given by \eqref{KKT_1}--\eqref{KKT_4}.
The chosen tracking type cost functional models the goal
of approximating
a desired fracture pattern by finding a suitable control.
However, the results extend to all Fréchet differentiable cost functionals.
We start by defining the required function spaces and cones.
Then, in \cref{Regularity of the upper level NLP},
we characterize the regular points for the upper level NLP
and conclude with the full KKT system.

Note that we drop the complementarity condition \eqref{KKT_5}
in order to obtain a Banach space NLP rather than a Banach space MPCC
(mathematical program with complementarity constraints) on the upper level.
This optimal control NLP is already novel and hard to solve.
The corresponding MPCC would not admit any regular point and would be
significantly more complicated both in theory and computation.
Of course, after solving the NLP it is easily checked
whether the omitted complementarity condition holds anyway,
i.e., whether we have a physically valid solution.

\subsection{Problem statement}

The upper level NLP is defined on the Hilbert space
$\mcY \define \CIQ \times Y \times  Z^*$
which is equipped with the natural norm
\begin{align*}
  \norm[\mcY]{(q,\boldu,\boldl)}^2
  &= \norm[W]{q}^2 + \norm[Y]{\boldu}^2 + \norm[Z^*]{\boldl}^2 \\
  &= \norm[L^2(I,Q)]{q}^2 + \norm[\Yu]{u}^2 + \norm[\Yphi]{\varphi}^2
    + \norm[\Vphi^*]{l_1}^2 + \norm[L^2(I,\Vphi^*)]{l_2}^2.
\end{align*}

To simplify the notation, we introduce a semilinear map
$a\colon \mcY \to Y^*$ representing the stationarity condition \eqref{KKT_4},
\begin{align*}
  a(q, \boldu, \boldl) \define {}
  &\prodIO{\gkappa(\varphi) \C e(u)}{e(P_u(\fcdot))}
    - \prodIN{q}{P_u(\fcdot)}
    + G_c \varepsilon \prodIO{\nabla\varphi}{\nabla P_\varphi(\fcdot)} \\
  &- \frac{G_c}{\varepsilon} \prodIO{1 - \varphi}{P_\varphi(\fcdot)}
    + (1 - \kappa) \prodIO{\varphi \C e(u)P_\varphi(\fcdot)}{e(u)} \\
  &- \dualprod{l_1}{P_\varphi(\fcdot)(0)}
    + \Idualprod{l_2(t)}{\pd_t P_\varphi(\fcdot)(t)}
    .
\end{align*}
The map $a$ is nonlinear in $\boldu$ and linear in $q$ and $\boldl$.
The cost functional is of tracking type, including a Tikhonov regularization term
which is beneficial for numerical stabilization \cite{KhiSteiWi22_JCP,NeiWiWo17,NeiWiWo19}.
For some desired spatial phase-field $\varphi_d \in \Vphi$, a nominal control $q_r \in Q$
and a Tikhonov parameter $\alpha >0$ the cost functional
$\J \colon W \times Y \to \R$ is defined as
\begin{equation}
\label{NLP_cost}
 \J(q, \boldu) \define \frac12 \int_I
    \bigl( \norm[\Omega]{\varphi(t) - \varphi_d}^2
    + \alpha \norm[\Gamma_N]{q(t) - q_r}^2 \bigr) \dt.
\end{equation}

\begin{Problem}
\label{prob:NLP}
Let the parameters $\varepsilon$, $\kappa$, $\mu$, $\lambda$ and $\varphi_0$ be as in
\cref{prob:NLP_E}. For given $\varphi_d \in \Vphi$, $q_r \in Q$ and $\alpha >0$,
find a control $q \in W$ and functions $\boldu \in Y$, $\boldl \in Z^*$
that minimize the cost functional $\J$ subject to \cref{KKT_1,KKT_2,KKT_3,KKT_4}:
\begin{equation}
  \label{NLP}
  \begin{aligned}
    \min_{(q, \boldu, \boldl) \in \mcY} \quad
    &\J(q, \boldu)  \\
    \stq
    & \kern0.7em \varphi(0) = \varphi_0 \text{ in }\Vphi, \\
    & {-}\dphi(t) \ge 0 \text{ in } \Vphi \aein I, \\
    & \kern-1.0em a(q, \boldu, \boldl) = 0 \text{ in } Y^*,\\
    & \kern0.7em l_2(t) \ge 0 \text{ in } \Vphi^* \aein I.
  \end{aligned}
\end{equation}
\end{Problem}
We define the constraints operator $\mcG \colon \mcY \to \mcZ$ as
\begin{align*}
  \mcG
  \define
  \begin{pmatrix}
    \mcG_{\setE} \\
    \mcG_{\setI}
  \end{pmatrix}, \quad
  \mcG_{\setE}(q, \boldu, \boldl)
  &\define
    \begin{pmatrix}
      \varphi(0) - \varphi_0 \\
      a(q, \boldu, \boldl)
    \end{pmatrix} \in \Vphi \times Y^*
  , \\
  \mcG_{\setI}(q, \boldu, \boldl)
  &\define
    \begin{pmatrix}
      -\dphi \\
      l_2
    \end{pmatrix}
  \in L^2(I,\Vphi) \times L^2(I,\Vphi^*)
  \qfor (q,\boldu,\boldl)\in \mcY.
\end{align*}
Consequently the image space of $\mcG$ reads as
$\mcZ \define \mcZ_1 \times \mcZ_2 \times \mcZ_3 \times \mcZ_4$ where
\begin{align*}
  \mcZ_1 &\define \Vphi,
  &\mcZ_2 &\define Y^*,
  &\mcZ_3 &\define L^2(I, \Vphi),
  &\mcZ_4 &\define L^2(I, \Vphi^*).
\end{align*}
The dual space is then
\begin{align*}
  \mcZ^*
  =
  \Vphi^* \times Y^{**} \times L^2(I, \Vphi^*) \times L^2(I, \Vphi^{**}  ),
\end{align*}
and the constraints cone becomes
$\mcK \define \mcK_1 \times \mcK_2 \times \mcK_3 \times \mcK_4$ where
\begin{align*}
  \mcK_1 &\define \set{0} \subset \Vphi,
  &\mcK_3 &\define \defset{v \in \mcZ_3}{v \ge 0} \subset L^2(I, \Vphi), \\
  \mcK_2 &\define \set{0} \subset Y^*,
  &\mcK_4 &\define \defset{v \in \mcZ_4}{v \ge 0} \subset L^2(I, \Vphi^*).
\end{align*}
Finally the dual cone $\mcK^*$ of $\mcK \subset \mcZ$ is defined as
\begin{align*}
 \mcK^*
  &\define \defset{\bs k = (k_1, k_2, k_3, k_4) \in \mcZ^*}
    {\bs k \boldz \ge 0 \tfor \boldz = (0, 0, z_3, z_4) \in \mcK} \\
  &~=\Vphi^* \times Y^{**} \times \biggl\{
    (k_3, k_4) \in L^2(I, \Vphi^*) \times L^2(I, \Vphi^{**})\colon \biggr. \\
  & \kern5em \biggl.
    \int_I \bigl( \dualprod{k_3}{z_3} + \bidualprod{k_4}{z_4} \bigr) \dt
    \ge 0  \tfor (z_3, z_4) \in \mcK_3 \times \mcK_4 \biggr\}.
\end{align*}

\subsection{Derivatives}

Next we compute the derivatives of
$\J\colon \mcY \to \R$ and $\mcG\colon \mcY \to \mcZ$
for all $(q,\boldu,\boldl)\in \mcY$ and
any given direction $h = (\delta q, \bolddeltau, \bolddeltal) \in \mcY$,
where $\bolddeltau = (\delta u, \delta\varphi)$
and $\bolddeltal = (\delta l_1, \delta l_2)$.
The quadratic functional $\J$ is infinitely often Fréchet differentiable as a functional
$\J \colon (\mcY,\norm[L^2]{\fcdot}) \to \R$ and hence also as a functional
$\J \colon (\mcY,\norm[\mcY]{\fcdot}) \to \R$:
\begin{align*}
  \J(q,\boldu,\boldl)
  &=
  \frac12 \left( \normIO{\varphi - \varphi_d}^2 +
  \alpha \norm[I \times \Gamma_N]{q - q_r}^2 \right),
  \\
  \J'(q,\boldu,\boldl) h
  &=
  \prodIO{\varphi - \varphi_d}{\delta\varphi} +
  \alpha \prodIN{q - q_r}{\delta q},
  \\
  \J''(q,\boldu,\boldl)(h_1, h_2)
  &=
  \prodIO{\delta\varphi_1}{\delta\varphi_2} +
  \alpha \prodIN{\delta q_1}{\delta q_2},
\end{align*}
and $\J^{(k)} = 0$ for $k > 2$.
For the constraint map $\mcG$ we obtain the derivative
\begin{equation*}
  \mcG'(\bq, \bbu, \bbl)(\delta q, \bolddeltau, \bolddeltal)
  =
  \begin{pmatrix}
    \delta\varphi(0) \\
    a'(\bq, \bbu, \bbl)(\delta q, \bolddeltau, \bolddeltal) \\
    -\delta\dphi \\
    \delta l_2
  \end{pmatrix}
  ,
\end{equation*}
and we show in the following proposition that this actually is a Fréchet derivative,
focusing on the only nontrivial component.

\begin{proposition}\label{prop:a'}
 The semilinear form $a\colon \mcY \to Y^*$ is Fréchet differentiable
 at every point $(q,\boldu,\boldl) \in \mcY$. For each direction
 $(\delta q, \bolddeltau, \bolddeltal) \in \mcY$, the derivative reads
\begin{align*}
  a'(q, \boldu, \bs l)(\delta q, \bolddeltau, \bolddeltal)
  &= \prodIO{\gkappa(\varphi) \C e(\delta u)}{e(P_u(\fcdot))} \\
  &\quad+ 2 (1-\kappa) \prodIO{\delta\varphi \varphi \C e(u)}{e(P_u(\fcdot))} \\
  &\quad- \prodIN{\delta q}{P_u(\fcdot)}
    + G_c \varepsilon \prodIO{\nabla\delta\varphi}{\nabla P_\varphi(\cdot)} \\
  &\quad+ \frac{G_c}{\varepsilon} \prodIO{\delta\varphi}{P_\varphi(\fcdot)}
    + 2(1-\kappa) \prodIO{\varphi P_\varphi(\fcdot) \C e(\delta u) }{e(u)} \\
  &\quad+ (1-\kappa) \prodIO{\delta\varphi P_\varphi(\fcdot) \C e(u) }{e(u)} \\
  &\quad- \dualprod{\delta l_1}{P_\varphi(\fcdot)(0)}
    + \Idualprod{\delta l_2}{\pd_t P_\varphi(\fcdot)}
      .
\end{align*}
\end{proposition}
\newcommand\Psiboldu[1][]{\Psi_{\boldu\optindex{#1}}}
\newcommand\Psiboldl[1][]{\Psi_{\boldl\optindex{#1}}}
\begin{proof}
We proceed as in the proof of \cref{prop:f'+f''}.
Straightforward computations confirm that the expression for $a'$ is indeed the
Gâteaux derivative of $a$.
For arbitrary $(q,\boldu,\boldl)$, $(\Psi_q,\Psi_{\boldu},\Psi_{\boldl})$, 
$(\delta q, \bolddeltau, \bolddeltal) \in\mcY$ and $\boldPhi \in Y$
we apply \cref{lem:estimatescalarproduct} and the identity
$\gkappa(\varphi+\Psiphi)=\gkappa(\varphi)+(1-κ)(2φ\Psiphi+\Psiphi^2)$
to obtain
\begin{align*}
  \abs{a'
  &(q + \Psi_q, \boldu + \Psiboldu, \bs l + \Psiboldl)
    (\delta q, \bolddeltau, \bolddeltal)(\boldPhi)
    - a'(q, \boldu, \bs l)(\delta q, \bolddeltau, \bolddeltal)(\boldPhi)} \\
  &= \bigl|
    (1-κ)\prodIO{(2φ\Psiphi+\Psiphi^2)\C e(\delta u)}{e(\Phiu)}
    + 2(1-κ)\prodIO{φ \delta φ \C e(\Psiu)}{e(\Phiu)} \\
  &\quad\;+ 2(1-κ)\prodIO{\Psiphi \delta φ \C e(u)}{e(\Phiu)}
    + 2(1-κ)\prodIO{\Psiphi \delta φ \C e(\Psiu)}{e(\Phiu)} \\
  &\quad\;+ 2(1-κ)\prodIO{φ\Phiphi \C e(\delta u)}{e(\Psiu)}
    + 2(1-κ)\prodIO{\Psiphi\Phiphi \C e(\delta u)}{e(u)}\\
  &\quad\;+ 2(1-κ)\prodIO{\Psiphi\Phiphi \C e(\delta u)}{e(\Psiu)}
    + 2(1-κ)\prodIO{\delta \varphi \Phiphi \C e(u)}{e(\Psiu)}\\
  &\quad\;+ (1-κ)\prodIO{\delta \varphi \Phiphi \C e(\Psiu)}{e(\Psiu)} \bigr| \\
  &\le (1-κ)C\bigl(
    \norm[\Linf(\IO)]{2φ\Psiphi+\Psiphi^2} \norm[Y]{\bolddeltau}
    + 2\norm[\Linf(\IO)]{φ \delta φ} \norm[Y]{\Psiboldu}\\
  &\kern5.5em+ 2 \norm[\Linf(\IO)]{\Psiphi \delta φ} \norm[Y]{\boldu}
    + 2 \norm[\Linf(\IO)]{\Psiphi \delta φ}
    \norm[Y]{\Psiboldu}\bigr) \norm[Y]{\boldPhi} \\
  &\quad+ (1-κ)C^2\bigl(
    2\norm[\Linf(\IO)]{φ} \norm[Y]{\bolddeltau}\norm[Y]{\Psiboldu}
    + 2\norm[\Linf(\IO)]{\Psiphi} \norm[Y]{\bolddeltau}\norm[Y]{\boldu}\\
  &\kern7em+ 2\norm[\Linf(\IO)]{\Psiphi}
    \norm[Y]{\bolddeltau}\norm[Y]{\Psiboldu}
    + 2\norm[\Linf(\IO)]{\delta φ} \norm[Y]{\boldu}\norm[Y]{\Psiboldu}\\
  &\kern7em+ \norm[\Linf(\IO)]{\delta φ}
    \norm[Y]{\Psiboldu}\norm[Y]{\Psiboldu}\bigr) \norm[Y]{\boldPhi},
\end{align*}
where $C$ is from \cref{lem:estimatescalarproduct}. Consequently it holds that
\begin{align*}
  \norm[Y^*]{a'
  &(q + \Psi_q, \boldu + \Psiboldu, \bs l + \Psiboldl)
    (\delta q, \bolddeltau, \bolddeltal)
    - a'(q, \boldu, \bs l)(\delta q, \bolddeltau, \bolddeltal)} \\
  &\le (1-κ)C\bigl(
    \norm[\Linf(\IO)]{2φ\Psiphi+\Psiphi^2} \norm[Y]{\bolddeltau}
    + 2C^2\norm[Y]{\boldu} \norm[Y]{\bolddeltau} \norm[Y]{\Psiboldu} \\
  &\kern5.5em+ 2C^2\norm[Y]{\Psiboldu}\norm[Y]{\bolddeltau} \norm[Y]{\boldu}
    + 2C^2 \norm[Y]{\Psiboldu}\norm[Y]{ \bolddeltau} \norm[Y]{\Psiboldu}\bigr) \\
  &\quad+(1-κ)C^2\bigl(
    2\norm[\Linf(\IO)]{φ} \norm[Y]{\bolddeltau}\norm[Y]{\Psiboldu}
    + 2\norm[\Linf(\IO)]{\Psiphi} \norm[Y]{\bolddeltau}\norm[Y]{\boldu}\\
  &\kern7em+ 2\norm[\Linf(\IO)]{\Psiphi}
    \norm[Y]{\bolddeltau}\norm[Y]{\Psiboldu}
    + 2C\norm[Y]{\bolddeltau} \norm[Y]{\boldu}\norm[Y]{\Psiboldu} \\
  &\kern7em+ C\norm[Y]{\bolddeltau}
    \norm[Y]{\Psiboldu}\norm[Y]{\Psiboldu}\bigr),
\end{align*}
and therefore  $a'(q + \Psi_q, \boldu + \Psiboldu, \bs l + \Psiboldl) \to
a'(q, \boldu, \bs l)$ in $L(\mcY,Y^*)$ as $(\Psi_q,\Psiboldu,\Psiboldl)\to 0$.
This implies continuity of the Gâteaux derivative, which ensures Fréchet differentiability.
\end{proof}
Given any multiplier $\bs \pi = (\pi_1, \pi_2, \pi_3, \pi_4) \in \mcK^*$,
we thus have
\begin{align*}
  \bs \pi \mcG'(\bq, \bbu, \bbl)(\delta q, \bolddeltau, \bolddeltal)
  &= \dualprod{\pi_1}{\delta \varphi(0)}
    + \sprod[Y^{**},Y^*]{\pi_2}{a'(\bq,\bbu,\bbl)(\delta q, \bolddeltau, \bolddeltal)} \\
  &\quad- \Idualprod{\pi_3}{\delta\dphi}
    - \int_I \bidualprod{\pi_4}{\delta l_2} \dt
    .
\end{align*}
\subsection{Regularity for the upper level NLP}
\label{Regularity of the upper level NLP}

A feasible point $(\bq, \bbu, \bbl) \in \mathcal M$ with
\begin{math}
  \mathcal M
  =
  \defset{(q, \boldu, \boldl) \in \mcY}{\mcG(q, \boldu, \boldl) \in \mcK}
\end{math}
will be regular for \cref{prob:NLP} if
\begin{equation*}
  \mcZ
  =
  \mcG'(\bq, \bbu, \bbl) \mcY - \mcK_{\mcG(\bq, \bbu, \bbl)},
\end{equation*}
where
\begin{math}
  \mcK_{\mcG(\bq, \bbu, \bbl)}
  =
  \mcK + \defset{\alpha (0, 0, -\dot\bphi, -\bl_2)}{\alpha \in \R}.
\end{math}
Thus, given $(\bq, \bbu, \bbl) \in \mathcal M$
and any $\boldz = (z_1, z_2, z_3, z_4) \in \mcZ$, we seek
$(\delta q, \bolddeltau, \bolddeltal) \in \mcY$,
$(k_3, k_4) \in \mcK_3 \times \mcK_4$,
and $\alpha \in \R$ such that
\begin{equation}
\label{eq:reg}
\begin{aligned}
  z_1 &= \delta\varphi(0) \text{ in } \Vphi, \\
  z_2 &= a'(\bq, \bbu, \bbl)(\delta q, \bolddeltau, \bolddeltal) \text{ in } Y^*, \\
  z_3 &= -\delta\dphi - k_3 - \alpha \dot\bphi \text{ in } L^2(I, \Vphi), \\
  z_4 &= -\delta l_2 - k_4 - \alpha \bl_2 \text{ in } L^2(I, \Vphi^{*}).
\end{aligned}
\end{equation}
The following proposition provides a sufficient condition
for the characterization of regular points.
\begin{proposition}
\label{prop:reg}
  Let $\bbu = (\bu, \bphi) \in Y$ be given such that
  \begin{equation}
    \label{eq:north}
    \C [\gkappa(\bphi)e(\wu) + 2(1-\kappa) \bphi e(\bu) \wphi]
    \not\perp \mathcal E
    \qfor \bs w \define (\wu, \wphi) \ne 0,
  \end{equation}
  where $\mathcal E \define \defset{e(v)}{v \in \Yu} \subset [L^2(I,L^2(\Omega))]^{2\times 2}$.
  Then the equation system \eqref{eq:reg} admits a solution
  $(\delta q, \delta u, \delta \varphi, \delta l_1, \delta l_2,k_3,k_4,\alpha) \in \mcY\times \mcK_3\times \mcK_4\times \R$.
\end{proposition}
\begin{proof} Let $\boldz = (z_1, z_2, z_3, z_4) \in \mcZ$ be arbitrary.
  Set $\delta\varphi(t) = z_1 - \int_0^t z_3(s) \ds$, $k_3 = 0$,
  $k_4 = 0$, $\delta l_2 = -z_4$, $\alpha = 0$, and $\delta q = 0$.
  Then the first, third and fourth equations of \eqref{eq:reg} hold and
  $\sup_{t\in I}\norm[\Vphi]{\delta\varphi(t)} < \infty$.
  It remains to show that we can find solution components
  $(\delta u, \delta l_1) \in Y_u \times \Vphi^*$
  for the second equation,
  \[
    a'(\bq, \bbu, \bbl)(\delta q, \bolddeltau, \bolddeltal)(\bs w) =
    \dualprod{z_2}{\bs w} \qfor \bs w \define (\wu, \wphi) \in Y.
  \]
  This is equivalent to solving $b((\delta u, \delta l_1),\bs w) = r(\bs w)$
  for all $\bs w \in Y$, where $b\colon \tilde\mcY  \times Y \to \R$
  with $\tilde\mcY \define \Yu \times \Vphi^*$
  denotes the bilinear form defined as
  \begin{align*}
    b((\delta u, \delta l_1), \bs w)
    &\define \prodIO{\gkappa(\bphi) \C e(\delta u)}{e(\wu)}
      + 2(1 - \kappa) \prodIO{\bphi \wphi \C e(\delta u)}{e(\bu)} \\
    &\kern1.5em- \dualprod{\delta l_1}{\wphi(0)}
      \qfor (δu,δl)\in \mcY, \ \bs w \in Y,
  \end{align*}
  and $r(\bs w)$ denotes the corresponding right hand side
  \begin{align*}
    r(\bs w) \define \dualprod{z_2}{\bs w}
    &- \Bigl(
      2 (1 - \kappa) \prodIO{\delta\varphi \bphi \C e(\bu)}{e(\wu)}
      + G_c \varepsilon \prodIO{\nabla\delta\varphi}{\nabla \wphi} \\[-\jot]
    &\kern2em + (1 - \kappa) \prodIO{\delta\varphi \wphi\C e(\bu)}{e(\bu)} \\
    &\kern2em + \frac{G_c}{\varepsilon} \prodIO{\delta\varphi}{\wphi}
      + \Idualprod{\delta l_2}{\partial_t \wphi} \Bigr),
      \quad \bs w \in Y.
  \end{align*}
  Clearly, $r\in Y^*$.
  In order to apply the Babuška-Lax-Milgram theorem \cite[Theorem 5.1.2]{Quarteroni1994},
  we have to show that $b$ is continuous and weakly coercive:
  \begin{align}
    \label{eq:coerA}
    \exists C > 0 \colon
    \sup_{\norm[Y]{\bs w}=1} \abs{b((\delta u, \delta l_1),\bs w)}
    &\ge C \norm[\tilde\mcY]{(\delta u, \delta l_1)}
      \qfor (δu,δl_1) \in \tilde\mcY, \\
    \label{eq:coerB}
    \sup_{\norm[\tilde\mcY]{(\delta u, \delta l_1)} = 1}
    \abs{b((\delta u, \delta l_1),\bs w)}
    &> 0 \qfor 0 \ne \bs w \in Y,
  \end{align}
  where $\norm[\tilde\mcY]{(\delta u,\delta l_1)}
  \define \norm{\delta u}_{Y_u} + \norm{\delta l_1}_{\Vphi^*}$.
  \newcommand\Lphi{U}%
  Setting $\Lphi \define \Linf(\IO)$ for brevity,
    we obtain continuity from the fact that there is $c>0$ such that
  \begin{align*}
    \abs{b((\delta u, \delta l_1), \bs w)}
    &\le \abs{\prodIO{\gkappa(\bphi) \C e(\delta u)}{e(\wu)}}
      + \abs{2(1 - \kappa) \prodIO{\bphi \wphi \C e(\delta u)}{e(\bu)}} \\
    &\quad+ \abs{\dualprod{\delta l_1}{\wphi(0)}} \\
    &\le c \bigl(
      \norm[\Lphi]{\bphi}^2 \normIO{∇\delta u}\normIO{∇\wu}
      + \norm[\Lphi]{\bphi} \normIO{∇\delta u}
      \normIO{∇\bu} \norm[\Lphi]{\wphi} \\
    &\qquad+ \norm{\delta l_1}_{\Vphi^*} \norm{\wphi(0)}_{\Vphi} \bigr)
  \end{align*}
  for all $(δu,δl_1,\bs w)\in \tilde\mcY\times Y$.
  With $c_1 \define \norm[\Lphi]{\bphi}$, $c_2 \define \normIO{∇\bu} $ 
  and $\tilde c=c(c_1^2+c_1c_2+1)$ it holds that
  \begin{align*}
    \abs{b((\delta u, \delta l_1), \bs w)}
    &\le c \bigl( c_1^2 \normIO{∇\delta u}\normIO{∇\wu}
      + c_1 c_2 \normIO{∇\delta u} \norm[\Lphi]{\wphi} \\
    &\kern2em + \norm[\Vphi^*]{\delta l_1} \norm[\Vphi]{\wphi(0)} \bigr) \\
    &\le c \bigl(
      c_1^2 \norm[\tilde\mcY]{(\delta u, \delta l_1)} \norm[Y]{\bs w}
      + c_1 c_2 \norm[\tilde\mcY]{(\delta u, \delta l_1)} \norm[Y]{\bs w} \\
    &\kern2em
      + \norm[\tilde\mcY]{(\delta u, \delta l_1)} \norm[Y]{\bs w} \bigr) \\
    &\le \tilde{c} \norm[\tilde\mcY]{(\delta u, \delta l_1)} \norm[Y]{\bs w}.
  \end{align*}
  Now the first coercivity condition \cref{eq:coerA} can be seen to hold by setting
  $(\wu,\wphi) = (\delta u / \norm[\Yu]{\delta u}, 0)$ in the supremum,
  \begin{align*}
    \sup_{\norm[Y]{(\wu,\wphi)}=1}
    &\bigl| \prodIO{\gkappa(\bphi) \C e(\delta u)}{e(\wu)}
      + (1 - \kappa) \prodIO{\bphi \wphi \C e(\delta u)}{e(\bu)} \\[-\jot]
    &+ \dualprod{\delta l_1}{\wphi(0)} \bigr| \\
    \ge
    & \Abs{ \frac{1}{\norm[\Yu]{\delta u}}
      \prodIO{\gkappa(\bphi) \C e(\delta u)}{e(\delta u)} }
      \ge \frac{C}{\norm[\Yu]{\delta u}} \norm[\Yu]{\delta u}^2=C\norm[\Yu]{\delta u},
  \end{align*}
  where $C$ in the last inequality is taken from Korn's inequality.
  To prove the second coercivity condition \cref{eq:coerB},
  we distinguish the cases $\wphi(0) = 0$ and $\wphi(0) \ne 0$.
  If $\wphi(0) \ne 0$, then there exists
  $\delta l_1 \in Z_1^*$ with $\norm[Z_1^*]{δl_1} = 1$
  such that \[\dualprod{\delta l_1}{\wphi(0)} \ne 0.\]
  By setting $\delta u = 0$ we have
  $\norm[\tilde\mcY]{(\delta u, \delta l_1)} = 1$,
  and it holds that
  \[
    \abs{b((\delta u, \delta l_1),\bs w)} =
    \abs{\dualprod{\delta l_1}{\wphi(0)}} > 0.
  \]
  If $\wphi(0) = 0$, we set $\delta l_1 = 0$ and conclude from \cref{eq:north} that
  there exists $\delta u \in \Yu$ with $\norm[\Yu]{\delta u} = 1$
  such that
  \begin{align*}
    \prodIO{e(\delta u)}
    {\C[\gkappa (\bphi) e(\wu) + 2(1 - \kappa) \bphi e(\bu) \wphi]}
    &\ne 0
  \end{align*}
  and thus
  \begin{align*}
  \abs{b((\delta u, \delta l_1),\bs w)} = \abs{\prodIO{\gkappa(\bphi) \C e(\delta u)}{e(\wu)}
    + 2(1 - \kappa) \prodIO{\bphi \wphi \C e(\delta u)}{e(\bu)}}>0.
  \end{align*}
  Consequently \cref{eq:coerA} and \cref{eq:coerB} hold and,
  due to the Babuška-Lax-Milgram theorem,
  there exists a solution $(\delta q, \bolddeltau,\bolddeltal) \in \mcY$
  of the second equation of \cref{eq:reg},
  and hence a solution of the entire system \cref{eq:reg}.
\end{proof}
By \cref{prop:reg} we can ensure that every feasible point $(\bq, \bbu,\bbl) \in \mathcal M$ which
fulfills condition \cref{eq:north} is regular for \cref{prob:NLP}.

\begin{Remark}
  It is apparent that every pair $(\bu,\bphi)$
  with $e(\bu) = 0$ or $\bphi = 0$ violates the condition \cref{eq:north}:
  simply choose $\bs w = (0,\wphi)$ with $\wphi \ne 0$.
  From a mechanics viewpoint, we first notice
  that a similar condition for the phase-field part exists,
  in which $2\C (1 - \kappa) \bphi e(\bu) \wphi$ is interpreted
  as some force that drives the fracture field
  (\cite[Section~3]{MieWelHof10b}; see also \cite[Section 4.5.3]{Wi20_book}).
  This interpretation is related to the complementarity condition that
  relates the bulk energy to crack growth.
  The situation is similar in \cref{eq:north}
  except that $\C (1 - \kappa) \bphi e(\bu) \wphi$ acts as
  a right hand side driving force for the displacement equation.
  This can be compared with classical elasticity, where the conservation of momentum
  is driven by right hand side volume and traction forces.
  Clearly, when the right hand sides are zero, and having
  $u\in V_u$, i.e., $u=0$ on $\Gamma_D$,
  we obtain trivial solutions,
  which are not of interest from the mechanics viewpoint.
  Thus, condition \cref{eq:north} simply
  excludes mechanically irrelevant solutions.
\end{Remark}
\subsection{Optimality conditions}
In order to formulate optimality conditions, we finally translate \cref{KKT} to the setting of problem \eqref{NLP}.
\begin{proposition}
Let $(\bq,\bbu,\bbl)\in \mathcal M$ be a regular minimizer of \eqref{NLP}.
Then there is a multiplier $\bs \pi = (\pi_1, \pi_2, \pi_3, \pi_4) \in \mcZ^*$ such that the KKT system
  \begin{align}
    \label{KKTN_1}\tag{KKTN\,1}
    \bphi(0) &= \varphi_0 \text{ in }  \Vphi, \\
    \label{KKTN_2}\tag{KKTN\,2}
    -\dot\bphi(t) &\ge 0 \text{ in } \Vphi \aein I, \\[\jot]
    \label{KKTN_3}\tag{KKTN\,3}
    a(\bq, \bbu, \bbl) &= 0 \text{ in }  Y^*, \\
    \label{KKTN_4}\tag{KKTN\,4}
    -\bl_2(t) &\ge 0 \text{ in } \Vphi^{*} \aein I, \\[\jot]
    \label{KKTN_5}\tag{KKTN\,5}
    \pi_3(t) & \ge 0 \text{ in }  \Vphi^{*} \aein I, \\
    \label{KKTN_6}\tag{KKTN\,6}
    \pi_4(t) & \ge 0 \text{ in }  \Vphi^{**} \aein I,
  \end{align}
  \vspace{-4ex}
  \begin{align}
    \notag
    \int_I \bigl[ \prodO{P_\varphi(\fcdot)}{\bphi - \varphi_d}
    + \prodO{P_q(\fcdot)}{\bq - q_r}\bigr] \dt \bigr. \kern6.5em  \\
    \notag {}
    - \dualprod{\pi_1}{P_\varphi(\fcdot)(0)}
    -\sprod[Y^{**},Y^*]{\pi_2}{a'(\bq,\bbu,\bbl)(P_q(\fcdot), P_{\boldu}(\fcdot), P_{\boldl}(\fcdot))} \kern1.5em \\
    \label{KKTN_7}\tag{KKTN\,7} \bigl. {}
    + \int_I \bigl( \dualprod{\pi_3}{\pd_t P_\varphi(\fcdot)}
    + \bidualprod{\pi_4}{P_{l_2}(\fcdot)} \bigr) \dt &= 0, \\
    \notag {}
      \dualprod{\pi_1}{\bphi(0) - \varphi_0}
    + \sprod[Y^{**},Y^*]{\pi_2}{a(\bq, \bbu, \bbl)} \kern1.5em \\
    \label{KKTN_8}\tag{KKTN\,8} \bigl. {}
    - \int_I \bigl( \dualprod{\pi_3}{\dot \bphi}
    - \bidualprod{\pi_4}{\bl_2} \bigr) \dt &=0
  \end{align}
  is satisfied, where we have used the projections defined by
  $P_q(\boldPhi) = \Phi_q$,
  $P_{\boldu}(\boldPhi) = \Phi_{\boldu}$,
  $P_{\boldl}(\boldPhi) = \Phi_{\boldl}$,
  $P_{l_2}(\boldPhi) = \Phi_{l_2}$.
\end{proposition}
 \begin{proof}
  Conditions \eqref{KKTN_1}--\eqref{KKTN_4} are just feasibility
  $(\bq,\bbu,\bbl)\in \mathcal M$ for the primal variables.
  Conditions \eqref{KKTN_5} and \eqref{KKTN_6} are feasibility
  for the multiplier, that is, $\bs \pi \in \mcK^*$.
  The stationarity condition
  $\J '(\bq, \bbu, \bbl) - \bs \pi \mcG'(\bq, \bbu, \bbl) = 0 \in \mcY^*$
  from \cref{KKT} is presented in \eqref{KKTN_7}.
  At last \eqref{KKTN_8} is the complementarity condition
  $\bs \pi \mcG(\bq, \bbu, \bbl) = 0$, again asserted by \cref{KKT}.
\end{proof}
\begin{Remark}
 For the sake of clarity, throughout the article,
 we have kept $\Vphi^*$ and $\Vphi$ separate,
 although, being a Hilbert space and its dual,
 they are isomorphic and could be identified with each other.
\end{Remark}

\section{Conclusions}
In this paper, we rigorously established a space-time
phase-field fracture complementarity model in
combination with an optimal control problem.
By formulating phase-field fracture as an abstract NLP in Banach spaces,
a complementarity system was obtained.
This derivation includes all cones necessary to characterize the multiplier.
Within this formulation the crack irreversibility was treated
as an inequality constraint for the
time derivative of the phase-field.
Hence a careful selection of a suitable functional framework was necessary
to obtain regularity results for the lower-level phase-field NLP.
In this process, all required differentiability results,
i.e., Fréchet differentiability of the energy and the constraints, were rigorously
shown.
In \cref{sec_further_opt_cond}, we discussed necessary optimality
conditions of second order together with first and second order sufficient conditions.
The KKT system resulting from the lower level problem then served as the constraint
for the optimal control problem, which is again formulated as an abstract NLP in Banach spaces.
Under certain conditions, regularity results for the higher level NLP were shown.
These conditions were then brief\/ly interpreted from a mechanical viewpoint.
Finally, we presented the full KKT system for the optimal control problem.

\section*{Acknowledgements}
The first and fourth author are partially funded by the
      Deutsche Forschungsgemeinschaft (DFG, German Research Foundation)
      Priority Program 1962 (DFG SPP 1962) within the subproject
      \emph{Optimizing Fracture Propagation using a Phase-Field Approach}
      with the project number 314067056.
      The third author is funded by the DFG -- SFB1463 -- 434502799.


\bibliographystyle{abbrv}

\end{document}